\documentclass[10pt,a4paper,twoside,leqno]{bull}
\usepackage{amsfonts,amsmath}
\usepackage{amssymb}

\usepackage{graphicx}

\font\Bigtit=cmr10 scaled \magstep 4
\font\ebf=cmbx8

\begin{document}

\thispagestyle{empty}

\begin{flushright}
PL ISSN 0459-6854
\end{flushright}
\vspace{0.5cm}
\centerline{\Bigtit B U L L E T I N}
\vspace{0.5cm}
\centerline{DE \ \  LA \ \  SOCI\'ET\'E \ \  DES \ \  SCIENCES \ \ ET \ \ DES \
\ \ LETTRES \ \ DE \ \ \L \'OD\'Z}
\vspace{0.3cm}
\noindent 2009\hfill Vol. LIX
\vspace{0.3cm}
\hrule
\vspace{5pt}
\noindent Recherches sur les d\'eformations \hfill Vol. LIX, no.\,1
\vspace{5pt}
\hrule
\vspace{0.3cm}
\noindent pp.~103--116



\vspace{0.7cm}

\noindent {\it Andrzej Krzysztof Kwa\'sniewski}

\vspace{0.5cm}

\noindent {\bf  ON COMPOSITIONS OF NUMBERS AND GRAPHS
}

\vspace{0.5cm}

\noindent {\ebf Summary}

{\small The main purpose of this note is to pose a couple of
problems which are easily  formulated  thought some seem to be not yet
solved. These problems are of general interest for discrete mathematics
including  a new twig of a bough of theory of graphs i.e. related with
given graph compositions. The problems result from and are served in the
entourage of series of exercises with hints based predominantly on
Knopfmacher et all. recent papers.

 }

\vspace{0.6cm}

\section {Number compositions}

\begin{quote}
{\it ``The number of compositions of $n$ into $k$ parts -- all
\textbf{distinct = ?} ... that is the question ...'' }
\end{quote}

Let us start by recalling indispensable notions thus establishing notation
and terminology. Following  Clark Kimberling  [1]  we define natural
number compositions and their family cardinality $C(n,k,a,b)$ as follows.

\vspace{2mm}
\noindent \textit{Definition 1.1}
The composition of the natural number $n$ is the vector
$$\left\langle i_1,i_2,...,i_k \right\rangle$$

\noindent solution of  the
Diophantine equation  $$ (D^*) \quad \quad \quad        i_1 + i_2
+...,i_k = n,  \ \ \  a \leq i_1 , i_2,...,i_k \leq b. $$

\vspace{0.2cm}
\noindent $C(n,k,a,b)$ = the number of vector solutions of $(D^*)$

\vspace{0.2cm}
\noindent \textit{Exercise 1.1.}  Consider $C(n,k,0,\infty)$ and
$C(n,k,1,\infty)$. Show that
$$
C(n,k,0,\infty) = \left(\begin{array}{c} n+k-1 \\k-1 \\
\end{array}\right), \qquad
C(n,k,1,\infty) = \left(\begin{array}{c} n-1 \\ k-1 \\
\end{array}\right)$$

\noindent \textbf{Hint.}  As for $C(n,k,0,\infty)$  consider below the obvious
Manhattan paths coding
by vector solutions of $(D*)$. Naturally $C(n,k,0,\infty)$ = number of
Manhattan paths from $A$ to $B$ covered according to the rule:

$n$ steps \textbf{right} $\rightarrow$ and $k-1$ steps \textbf{up} $\uparrow$
\begin{center}
(there are $k$  levels or $k$ floors or $k$ storeys)\\
with coding: $i_s$ = number of steps right at the $s$-th level (floor)
\end{center}

\vspace{2mm}

\begin{figure}[ht]
\begin{center}

\includegraphics[width=75mm]{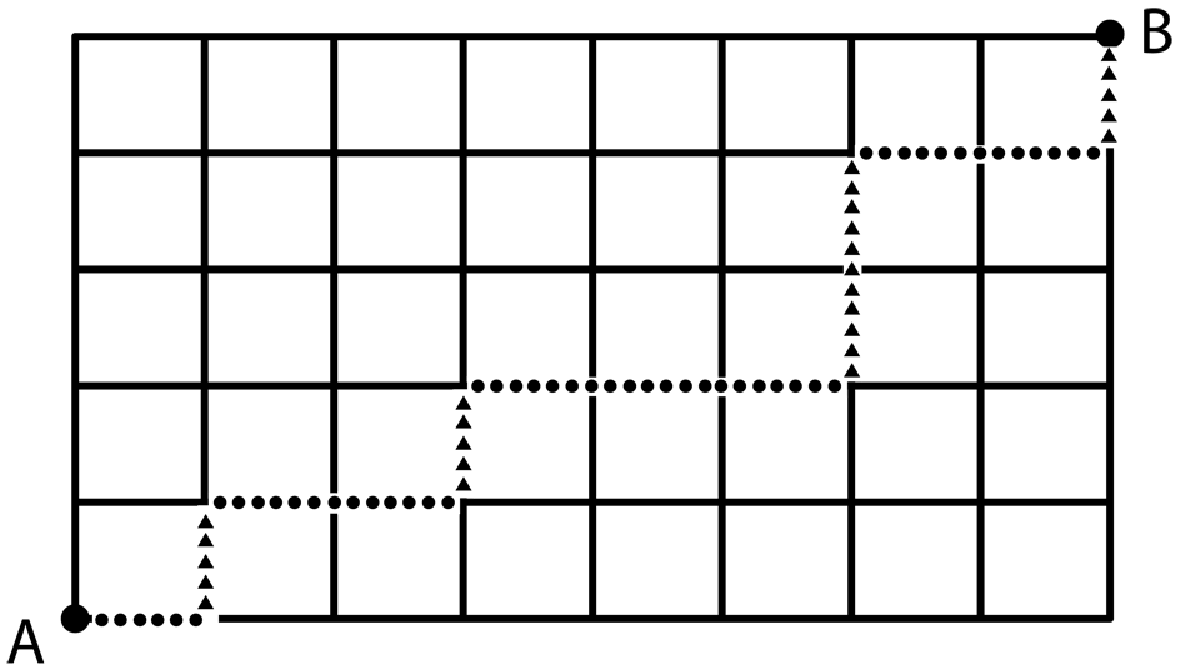}
\caption{ The path code $\left\langle i_1,i_2,...,i_6\right\rangle
=\left\langle 1,2,3,0,2,0\right\rangle.$ \label{fig:1}}
\end{center}
\end{figure}

Now let us go from $A$ to $B$ via a shortest way -- their number is the
number of coding strings i.e.
$$ \left(\begin{array}{c}
n+k-1 \\
n   \\
\end{array}\right) = \left(\begin{array}{c}
n+k-1 \\
k-1   \\
\end{array}\right). $$

\noindent \textbf{Of course:}
\begin{quote}
\emph{Compositions with the restriction $i_1\geq i_2\geq...\geq i_k \geq a$ are partitions.
\textbf{Note} that parts of compositions might be zero if $a = 0$.}
\end{quote}

\noindent Consider now partitions and compositions with pairwise distinct \emph{nonzero parts}.

\vspace{2mm}

\noindent \textit{Definition 1.2} Let us introduce the notation: \\
$\Pi[n,k]$ = the number of  \textbf{partitions} of $n$ into $k$ \textbf{distinct nonzero parts}. \\
$C[n,k]$ = the number of \textbf{compositionss} of $n$ into  $k$ \textbf{distinct nonzero parts}.

\vspace{4mm}

\noindent \textit{Exercise 1.2.}
\noindent Prove that $\Pi [0,0]= 1,\  \Pi[n,k] = 0$ if $n < 0$ and
$$\Pi [n,k]=\Pi [n-k,k]+ \Pi [n-k,k-1].$$

\vspace{2mm}

\noindent \textbf{Solution.}  \textbf{Start.} Indeed.\\
$\Pi[n-k,k-1]$ = the number of \textbf{partitions} of $n-k$ into $k-1$ distinct
non-zero parts = the number of \textbf{partitions} of $n$ into $k$  distinct
non-zero parts with the \textbf{smallest part equal to one}.
Indeed -- just look at this:
$$ 5 + 4 + 3 + 2 + 1 \quad\Leftrightarrow\quad  (5-1) + (4-1)  + (3-1) + (2-1) \  [+ \ zero]$$

\begin{figure}
\begin{center}
	\includegraphics[width=85mm]{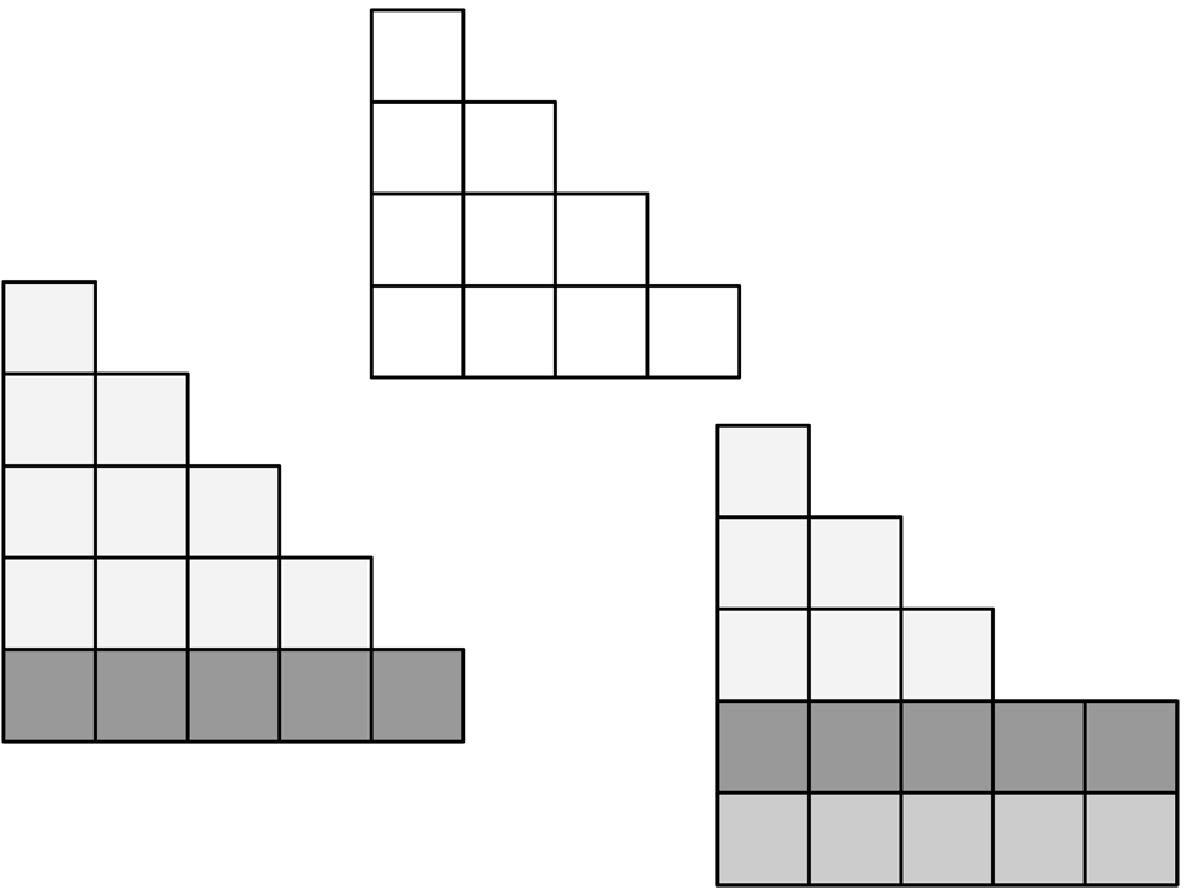}

\vspace{3mm}

	\caption{description \label{fig:2}}
\end{center}
\end{figure}

\noindent $k$ parts of  $n$ with smallest part equal to one $\Leftrightarrow k-1$ parts of $n-k$.
Here equivalence $\Leftrightarrow $ is achieved via cutoff of $k$ ones
(``unite boxes'') -- see Fig.\ref{fig:2}.

\vspace{2mm}

\vspace{1mm}

\noindent Now cut off $k$ unite boxes from $k$ parts with smallest part \textbf{greater than one}.
Then
$\Pi [n-k,k]$ = the number of \textbf{partitions} of $n-k$ into $k$ distinct non-zero parts
= the number of partitions of $n$ into $k$ distinct non-zero parts with the smallest part
\textbf{greater than one}.
Indeed -- just look at this:
$$ 6 + 4 + 3 + 2  \quad\Leftrightarrow\quad  (6-1) + ( 4-1)  + (3-1) + (2 - 1 ),$$

\noindent so $k$ parts of  $n$ with smallest part greater than one $\Leftrightarrow k$
distinct parts of $n-k.$

\vspace{2mm}

\noindent \textit{Exercise 1.3.}
\noindent Print a finite part of the Pascal -- like triangle given by the infinite array
with $\Pi [n,k]$ as matrix elements. (See: Mathemagics in \\ $http://ii.uwb.edu.pl/akk/index.html$).\\
Of course:
\begin{center}
the number of compositions $C[n,k] = k!$ number  of partitions  $\Pi [n,k]$.
\end{center}
Of course (Why \emph{"Of course"}?) $C[0,0] = 1$, $C[n,k] = 0$ if $n < 0$,
\begin{displaymath}
	C[n,k] = C[n-k,k] + k C[n-k,k-1].
\end{displaymath}
Why ``Of course''? Because there are $\uparrow$ $k$ places where from might
be found one box to be cut off. \\
\textbf{Note:} $C[n-k,k-1]$ = the number
of \textbf{compositions} of $n$ into $k$ distinct non-zero parts with the
smallest part \textbf{equal to one}.\\
$C[n-k,k]$ = the number of
\textbf{compositions} of $n$ into $k$ distinct non-zero parts with the
smallest part \textbf{greater than one}.

\vspace{2mm}

\noindent \textit{Exercise 1.4.}
\noindent Print a finite part of the Pascal-like triangle given by the
infinite array with $C[n,k]$
as matrix elements.
\vspace{2mm}

\begin{center}
\noindent \textbf{Compositions of $n$ into $k$ distinct parts problem I.}\\
\vspace{1mm}
\noindent The number of compositions of $n$ into $k$ distinct parts = ?\\
\emph{Find any compact formula for the answer.}
\end{center}

\emph{Any compact formula} as for example we have for Stirling numbers,
Newton binomial or Gauss $q$-binomial numbers etc. (See: Mathemagics in\\
\emph{http://ii.uwb.edu.pl/akk/}).

\section{Graph composition i.e.  ``comppartition'' of a graph  --\\ Basic
enumerations}

In this section proofs and hints follow [2].
In order to establish notation let us recall the meaning of the basic notions to be
used in what follows.   We shall consider finite, undirected, labeled graphs, with
no loops or multiple edges.
The edge with endpoints  $v_1$ and $v_2$ is  $\langle v_1, v_2\rangle$. The set of
vertices of graph $G$ is denoted by and the set of edges is denoted  by $E(G)$.
The notion we need now is that of an \textbf{induced} subgraph of the given graph G.

A subgraph $H$ of a graph $G$ is said to be induced if, for any pair of vertices $x$
and $y$ of $H$, $\langle v_1, v_2\rangle$ is an edge of $H$ if and only if
$\langle v_1, v_2\rangle$ is an edge of $G$. In other words, $H$ is an
\textbf{induced} subgraph of $G$ if it has the most edges from $E(G)$ over the
same vertex set  $V(H)$.
The idea of specific vertex partition of a graph - named  by A.\,Knopfmacher and
M.\,E.\,Mays [2] -- the graph composition, generalizes both ordinary compositions
of positive integers and partitions of finite sets.  As for the name we have
experienced that it  might be nowadays somewhat misleading therefore let us make
aware of two notifications  right after the Definition\,2.1. and related observation.

\vspace{0.2cm}
\noindent \textit{Definition 2.1} ({``comppartition'' = composition in} [2]){
Let $G = \langle E(G), V(G)\rangle$ denotes a labelled graph.
A \textbf{comppartition} of $G$  is a \textbf{partition} of the vertex set $V(G)$
into vertex subsets of
the \textbf{connected induced} \emph{subgraphs} of $G$ i.e. such a vertex
comp-partition of $G$ provides a \textbf{set} of \textbf{connected included} subgraphs of
$G \rightarrow  {G_1, G_2, ...,G_m}, G_i = \langle E(G_i), V(G_i)\rangle, i = 1,...,m$.}

\vspace{2mm}

\noindent \textit{Observation:}
It is important that $G_i$ = $\langle E(G_i), V(G_i)\rangle$ are \textbf{induced}
subgraphs because the \textbf{same} vertex subset may be spanned by different
\textbf{edge} subsets - therefore to the same \textbf{comppartition} of $G$ which
is a partition of vertex set $V(G)$ there might correspond   different families
$\{G_1, G_2, ...,G_m\} \neq \{ G'_1, G'_2, ..., G'_m\}$ of connected subgraphs.

\begin{figure}[ht]
\begin{center}
	\includegraphics[width=65mm]{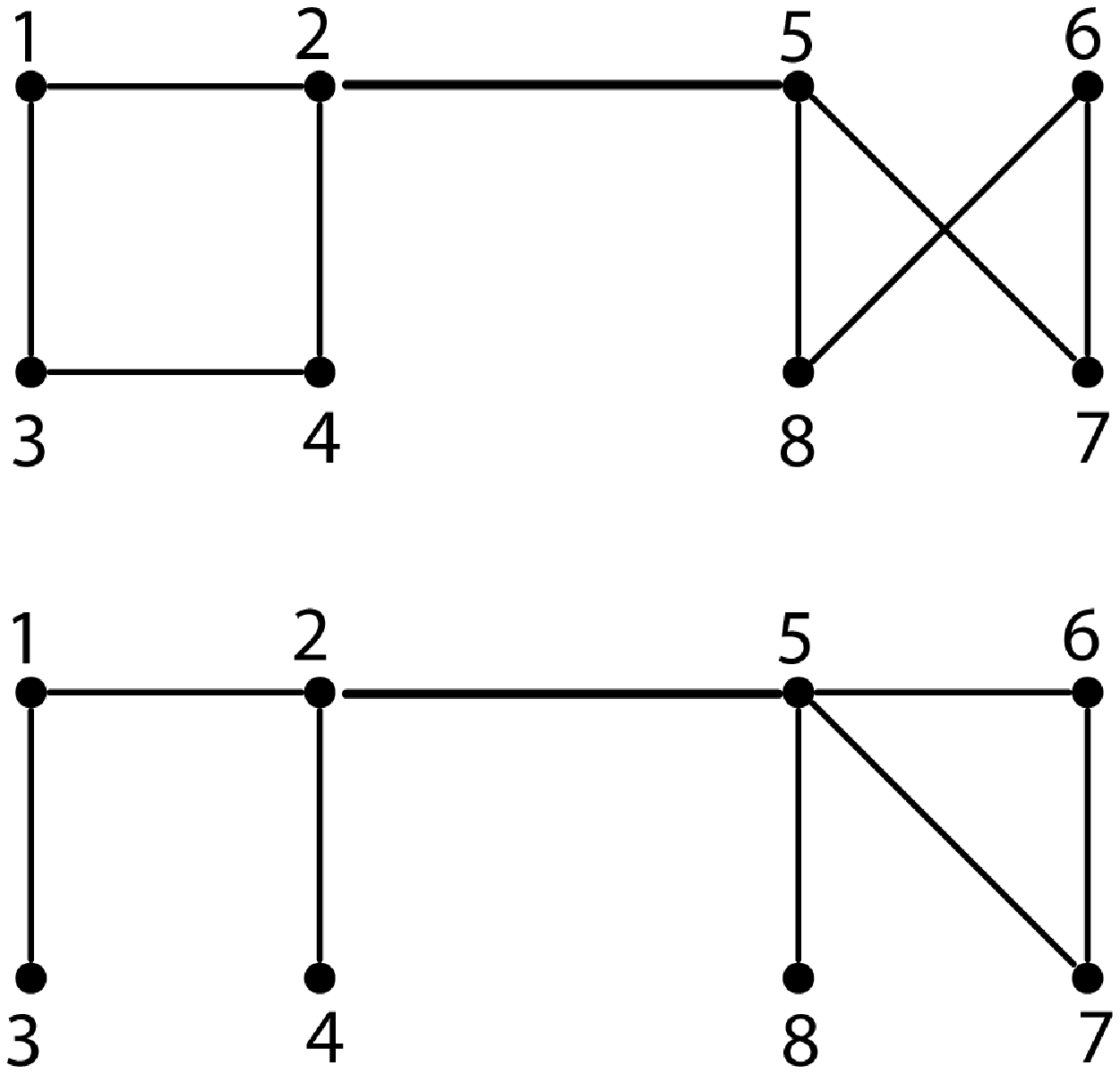}
	\caption{description \label{fig:3}}
\end{center}
\vspace{8mm}
\end{figure}

\begin{figure}[ht]
\begin{center}
	\includegraphics[width=75mm]{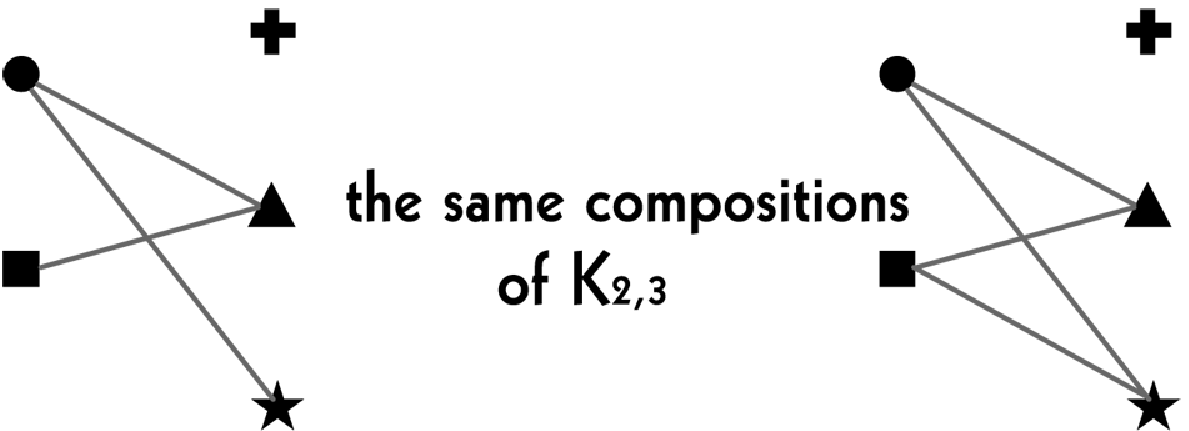}
	\caption{Bipartite complete  graph    $K_{2,3}$ compositions` example \label{fig:4}}
\end{center}
\end{figure}

\begin{quote}
$K_{n,m}$ complete biparite graph  has  $n\cdot m$  edges  linking  the first row $n$
dots in such a way that each  dot  of  this row is linked to every one of the second
row of $m$ dots (dots = vertices).
\end{quote}

\noindent \textbf{1. Note}
that the graph composition-partition introduced by A. Knopfmacher and M.\,E.\,Mays  is
not the textbook composition  defined as  follows. The composition $G = G_1[G_2]$ of
graphs $G_1$ and  $G_2$ with disjoint point sets $V_1$ and $V_2$ and edge sets $X_1$
and $X_2$ is the graph with  point vertex $V = V1\times V2$  and such that $u = (u_1,u_2)$
adjacent with $v = (v_1,v_2)$ whenever $[u_1$ adj $v_1]$ or  $[u_1 = v_1$ and $u_2$ adj $v_2]$.
It is also called the graph lexicographic product.

\vspace{0.2cm}

\noindent \textbf{2. Note} that $2^8 = 256$ different types of graph products may
be defined  of any given two graphs $G$  and $H$ is a new graph whose vertex set is
$V(G)\times V(H)$ and where, for any two vertices $(g,h)$ and $(g',h')$ in the product,
the adjacency of those two vertices is determined entirely by the adjacency (or equality, or non-adjacency) of $g$ and $g'$, and that of $h$ and $h'$.

\noindent [\emph{http://mathworld.wolfram.com/GraphComposition.html}] Neither of them is
``graph composition''  from [2] of course.

\vspace{0.3cm}

\noindent \textit{Notation:}
\noindent $C(G)$ = the number of distinct compositions of a graph $G$.

\begin{center}
 \textbf{Why the name: \textit{composition} of $G$ in [2] ?}\\
$\downarrow$
\end{center}

\noindent \textit{Exercise 2.1.}
Let $P_n$ be the path with $n$ vertices. Prove that $C(P_n)= 2^{n-1}, n>0,
C(P_0) \equiv 1$.
\begin{center}
 \textbf{Why not the name: \textit{partition} of $G$ in [2] ?}\\
$\downarrow$
\end{center}

\vspace{3mm}

\noindent \textit{Exercise  2.2.}
Let $K_n$ be the complete graph on $n$ vertices. Prove that
$C(K_n) = B_n, n>0, C(K_0)\equiv 1$. $B_n$ = Bell numbers = exponential numbers.\\

\noindent\textbf{Note again:} the same \textbf{partition} of $G$ which is a partition of vertex set $V(G)$  - may give rise to different families $\{G_1, G_2, ...,G_m\} \neq  \{G~_1, G~_2, ...,G~_m\}$ of connected subgraphs.
\begin{center}
\textbf{Why not the name: \emph{comppartition} of  G   in what follows ?} \\
$\downarrow$
\end{center}

\noindent \textit{Comment.} \emph{These two above examples  are two extreme cases:}\\
no connected graph $G$  with $n$ vertices can have fewer than $C(P_n)$ and no connected graph $G$ with   $n$ vertices can have  more than $C(K_n)$ \textbf{comppartitions}. i.e.
$2^{n-1} \leq C(F_n) \leq B_n$, where  $F_n$ is any connected graph with $n$ vertices.

\vspace{3mm}

\noindent \textit{Exercise 2.3.}
Let $G =  G_1 \cup G_2$ and there are no edges from vertices of $G_1$ to vertices of $G_2$.
Then $C(G)= C(G_1) C(G_2)$. The same holds for $G_1$ and $G_2$ having exactly one vertex in common.
Prove or rather see this.\\

\noindent \textit{Note:} one obtains \emph{comppartitions} of $G$ by pairing
\emph{comppartitions} of $G_1$ and $G_2$ in all possible ways.

\vspace{3mm}

\noindent \textit{Exercise 2.4.}
Let $G =  G_1 \cup G_2$ and there is exactly one edge from a vertex of $G_1$ to a vertex
of $G_2$ whose removal disconnects $G$. Then $C(G) = 2 C(G1)C(G2).$ Prove or rather see this.

\noindent \textit{Note:}  Let $e$ be the distinguished edge between vertices $v_j$ and $v_j.$
Any composition of $G$ can be obtained in exactly two ways: either $e$ is included
to supply the component $v_i$ in $G_1$ and the component $v_j$ in $G_2$ or not.
Thus the count from Exercise 2.3. is now doubled.

\vspace{3mm}

\noindent \textit{Exercise  2.5.}
Let $T_n$ be any tree with $n$ vertices. Prove that $C(T_n) =  2^{n-1},  n> 0$.

\vspace{3mm}

\noindent \textit{Proof:}   Use induction. Consider  $T_{n+1}$. Remove one edge.
This disconnects $T_{n+1}$ into two parts for which the formula holds. By the result
of  Exercise 2.4. the proof is accomplished.

\vspace{3mm}

\noindent \textit{Exercise  2.6.}
\noindent Let $K^-_n$  be the complete graph on $n$ vertices with \emph{one} edge removed.
Prove that $C(K^-_n) = B_n - B_{n-2}, n > 1.$

\vspace{3mm}
\noindent \textit{Proof:} Let $e$ be the \emph{deleted} edge between vertices $v_j$ and $v_j$. Its
deletion affects a composition counted by $C(K_n)  = B_n$ only when \emph{the component}
containing vertices $v_j$ and $v_j$ consists of \emph{exactly these two} vertices $v_j$ and $v_j$.
Otherwise there is a by-pass in $K_n$ connecting these two vertices $v_j$ and $v_j$.
Therefore from the number of composition counted by $C(K_n)  = B_n$ one must substract
those compositions for which one of the partition component is  ${v_j, v_j}$.
This restriction rules out  exactly   $C(K_n - 2)= B_{n-2}$ compositions of $K_n$.

\vspace{3mm}

\noindent \textit{Exercise  2.7.}
\noindent Let $C_n$ be the cycle graph with $n$ vertices. Prove that $C(C_n) = 2^n - n, n> 0$.

\vspace{3mm}

\noindent \textit{Proof:} Delete any  edge. The resulting graph becomes  $P_n$
with $C(P_n) = 2^{n-1}.$
Any composition of $P_n$ is also a composition of $C_n$. The deleted edge may be reinserted,
providing a new composition of $C_n$ - previously not counted \emph{unless} the composition of $P_n$
had been obtained by deleting from $P_n$ either no edge or exactly one edge.

In these cases, reinserting the original deleted edge gives the same composition of $C_n$:
namely the composition  consisting of all $n$ vertices. Therefore the total number
of distinct compositions of $C_n$ is equal to  $C(C_n) = 2 2^{n-1} - n  = 2^n  - n.$

\vspace{3mm}

\noindent \textit{Exercise  2.8.}  (Theorem 9 in [3])
\noindent Let $L_n$ be the ladder cycle graph with $2n$ vertices and $3n-2$ edges as $L_n$ is the
product of a path of lengt $n$ and a path of length $2$.Prove that $C(L_n)$ does satisfy
the following recurrence.
$$ (L_1) = 2, C(L_2) = 12,    C(Ln)   = 6 C(Ln-1)  + C(Ln-2) \ for \  n>2.$$

\noindent \textit{Proof:}  see: Theorem 9  in  [3].

\vspace{3mm}

\noindent \textit{Exercise  2.9. Number of Ladders} with $n$ rungs ``problem''.\\
\noindent $C(L_n)$ = ? Find the Binet-like formula for the answer.
[Jacques Binet (1786-1856)]. Contact:  $http://mathworld.wolfram.com/LadderGraph.html$.

\section{ Number compositions  with parts constrained\\ by the leading
summand}

In this section we follow  [3].  Note:  here now \textbf{again}:\\
Composition of a natural number is an ordered partition of  a natural number.\\
\emph{Ordered partition} of a natural number is a constrained composition of a natural number.\\
(...constrained?  yes !  constrained by...  answer this question)

\vspace{0.2cm}

\noindent \textbf{Of course:} Compositions constrained by  the requirement
$i_1\geq i_2\geq...\geq i_k\geq a$ are partitions. \\
\noindent \textbf{Note} \emph{that parts of compositions might be zero if  $a = 0$} [1].

\vspace{0.2cm}
\noindent \textit{Exercise  3.1.}
Consider the number  $f_n(k)$  of compositions  of a natural number $n$  into $k$
parts with the strictly largest part in the first position i.e.

\begin{displaymath}
	\mathbb{(>)} \ \ \ \ \ i_1 + i_2 + ... + i_k = n, \ \ \ 1 \leq i_1, i_2, ..., i_k,
	\ \ \ i_1 > i_k \ for \ k>1
\end{displaymath}


\noindent \emph{Observe-show that the following formula gives  the  ordinary generating
function for these numbers  with  $k$  fixed : }
$$
	F_k(z) = \sum_{n\geq 0}{f_n(k)z^n} = \frac{(1-z)z^k}{1 - 2z + z^k}
$$

\noindent \textit{Proof:}
$$k\geq 2 \quad \sum_{n\geq 0}{f_n(k)z^n}
= z^k\sum_{n\geq 1}{\Phi_n^{k-1}z^n} = \frac{z^k}{1-z-z^2-...-z^{k-1}},$$

\noindent where higher order Fibonacci  sequences are defined
$$\Phi_{n+k-1}^{(k-1)} = \sum_{i=0}^{k-2}{\Phi_{n+i}^{(k-1)}}  ...$$

\noindent (... what are the initial values?)

\vspace{2mm}

\noindent \textit{Exercise  3.2.}
{Consider the number  $f^*_n(k)$  of compositions  of a natural number $n$ into $k$
parts with the largest part in the first position } i.e.
\begin{displaymath}
	\mathbb{(\geq)} \ \ \ \ \ i_1 + i_2 + ... + i_k = n, \ \ \ 1 \leq i_1, i_2, ..., i_k,
	\ \ \ i_1 \geq i_k \ for \ k>1
\end{displaymath}

\vspace{2mm}

\noindent \emph{Observe-show that the following formula gives  the  ordinary generating  function
for these numbers with $k$ fixed: }
$$
	F_k^*(z) = \sum_{n\geq 0}{f_n^*(k)z^n} = \frac{(1-z)z^k}{1 - 2z + z^{k+1}}
$$

\noindent \textit{Proof:}
$$k\geq 2 \quad \sum_{n\geq 0}{f_n^*(k)z^n} =
z^k\sum_{n\geq 1}{\Phi_n^{k}z^n} = \frac{z^k}{1-z-z^2-...-z^{k}},$$

\noindent  where
higher order Fibonacci  sequences are defined
$$\Phi_{n+k}^{(k)} = \sum_{i=0}^{k-1}{\Phi_{n+i}^{(k)}}  ...$$

\noindent  (... what
are the initial values?)

\vspace{2mm}

\vspace{0.2cm}
\noindent \textit{Exercise  3.3.}
{Consider the two following sequences of numbers  $\langle f_n\rangle_{n\geq 0}$  and
$\langle f_n^*\rangle_{n\geq 0}$  and their ordinary generating functions}
$$
F(z) = \sum_{n\geq 0}F_n(z) = \sum_{n\geq 0}f_n z^n \ \ \ and \ \ \
F(z) = \sum_{n\geq 0}F_n^*(z) = \sum_{n\geq 0}f_n^* z^n
$$

\noindent \emph{Observe-show that for $n\geq 2 \ f_{n+1} = f^*_n$ as $zF^*(z) = F(z) - z$,
(...what are the initial values?)}

\vspace{0.2cm}
\noindent \textbf{Naturally}\\

\begin{center}
$f_n = \sum_{k\geq 0}f_n(k)$ = the number  of  all compositions  of a natural number
$n$ with the strictly largest part in the first position.
\end{center}

\vspace{0.2cm}
\noindent \textbf{Compositions} of $n$ with strictly largest part in the first position
\textbf{problem II.}

\begin{center}
	$f_n = ?$  Find any compact formula for the answer.
\end{center}

\noindent A possible way to solve the problem II  is to solve the recurrence equation
from the next exercise.

\vspace{0.2cm}
\noindent \textit{Exercise  3.4.} \\
Observe-show that
$$f_n(k) = 2f_{n-1}(k) - f_{n-k}(k) + \delta_{n,k} - \delta_{n,k+1},$$

\noindent (...what are the initial values?)

\noindent \textbf{Hint:} use
$$F_k(z) = \sum_{n\geq 0}f_n(k)z^n = \frac{(1 - z)z^k}{1-2z+z^k}.$$

\begin{center}
	\textbf{Compositions} of $n$ into $k$ parts with strictly largest part\\
 in the first position \textbf{problem  III}. \\
	$f_n(k) = \ ? \ \ $ Find any compact formula for the answer.
\end{center}

\section{Compositions with distinguished part}

In this section we follow  [4].

\vspace{0.2cm}
\noindent \textit{Exercise  4.1.}
Let $C_k(n)$ be the number of compositions of $n$ in which at least one $k$ occurs.
Prove then that
$$\sum_{n\geq 0}C_k(n)z^n = \frac{z}{1-2z} - \frac{z - z^k + z^{k+1}}{1 - 2z + z^k
- z^{k+1}},\quad  n\geq k \geq 1.$$

\vspace{0.2cm}
\noindent \textit{Proof:} Consider $C^*k(n)$ to be the  number of compositions of $n$ in
which no part equals to $k$.  Consider $C^*_k(n,m)$ to be  the  number of compositions of
$n$ into $m$ parts with no part equal to $k$.

\vspace{0.2cm}
\noindent \textbf{Note:}  by the product law of generating functions
$$
	\sum_{n\geq 0}C^*_k(n,m)z^n = (z + z^2 + ... + z^{k-1} + z^{k+1} + ...)^m
$$

\noindent and \textbf{sum} now the above over $m$ via $\sum_{m=1}^{n}$ and then

\vspace{0.2cm}

\noindent \textit{Observe} that
$$\frac{z - z^k + z^{k+1}}{1 - 2z + z^k - z^{k+1}}
\equiv \frac{q}{1-q}.\qquad q =  ?$$

\noindent \textbf{That is all} ... as of course
$$\sum_{n\geq 0}C(n)z^n = ... = \frac{z}{1-2z}.$$

\noindent ... and $q=?$ Answer:
$$q = \frac{z - z^k + z^{k+1}}{1-z}.$$

\vspace{0.2cm}
\noindent \textit{Exercise  4.2.}
$C^*k(n)$ to be the number of compositions of $n$ in which no part equals to $k$.  \\
Then
$$
	C^*_k(n) = 2C^*_k(n-1) - C^*_k(n-k) + C^*_k(n-k+1), \ \ \ n\geq k + 2\\
$$
$$
(...\ what\ are\ the\ initial\ values?)
$$

\vspace{0.2cm}

\noindent \textit{Proof:} see [4] for two proofs of the above. One of them is combinatorial.

\begin{center}
	\textbf{Compositions} of $n$ with no $k \in \mathbb{N}$ summand  \textbf{problem  IV}.\\
	$C^*_k(n)= $ ? Find any compact formula for the answer.
\end{center}

\section {Number compositions  with all part distinct}

For this section see [5]. See also related [6]. We come back to the very first question
posed at the start of  this article.

\vspace{0.2cm}

\noindent ``The number of compositions of $n$ into $k$  \textbf{all} parts
\textbf{distinct =  ?} \emph{ ... that is the question...}''
Recalling now the Definition 2.1. and using the result of Exercise 2.1. we proved what follows.

\vspace{0.2cm}
\noindent \textit{Recall 5.1.} Prove that $C[0,0] = 1$, $C[n,k] = 0$ if $n < 0$ and
$$
	C [n,k] = C[n-k,k] + k C[n-k,k-1].
$$

\noindent Let us now introduce also  the overall number of compositions into distinct parts.

\vspace{0.2cm}

\noindent \textit{Definition 5.1.} \emph{ $C[n] = \sum_{k\geq 1}C[n,k]$ = the number of
all compositions of  $n$  into \textbf{distinct} parts.   }

\vspace{0.2cm}
\noindent \textbf{Information} from [5]:  one may prove that

$$
	C[z] = \sum_{n\geq 1}C[n] = \sum_{k\geq 1}\frac{k! z^{{k+1 \choose 2}}}{(1-z)(1-z^2)
	...(1-z^k)}
$$

\noindent So what about our question now?

\begin{center}
	\textbf{Compositions} of  $n$ into   \textbf{distinct}   parts  \textbf{problem  V}.\\
The number of compositions of $n$   into  \textbf{distinct}   parts i.e.  $C[n]   =  ?$\\
Find any compact formula for the answer.
\end{center}

\noindent To this end and  for any case we add also the miscellaneous Appendix on formulas
with intrinsic reference to unavoidable numbers` compositions  therein.

\section{Appendix  ``Disce puer''}

Consider Stirling numbers [Knuth notation]  of the second $\Big\{ {n \atop k} \Big\}$
and of the first kind $\Big[ {n \atop k} \Big]$ via compositions of the natural number
$n$ \textbf{summation}. (See: Mathemagics in \emph{http://ii.uwb.edu.pl/akk/}).

\vspace{0.2cm}

\noindent \textit{Exercise 1.A.}
Recall [7,8,9], prove adn compare

$$
	\Bigg\{ {n \atop k} \Bigg\} = \frac{1}{k!}\sum_{{i_1+i_2+...+i_k = n} \atop
	{0<i_1,...,i_k\leq n}}{n \choose {i_1,i_2,...,i_k}}
$$

\noindent where

$$
	{n \choose {i_1,i_2,...,i_k}} = \frac{n!}{i_1!i_2!...i_k!}
$$

\noindent and

$$
	\mathbf{(D^*)} \ \ \ \ i_1 + i_2 + ... + i_k = n, \ \ 0 < i_1,i_2,...,i_k \leq n
$$

\noindent thereby $\langle i_1,i_2,...,i_k\rangle$ is a vector solution of $(D^*)$ i.e.
a composition of the natural number $n$.

\vspace{0.2cm}

\noindent Now compare

$$
	\Bigg\{ {n \atop k} \Bigg\} = \frac{n!}{k!} \sum_{{i_1+i_2+...+i_k = n} \atop
	{0<i_1,...,i_k\leq n}} \frac{1}{i_1!i_2!...i_k!}
$$

\noindent with [8]

$$
	\Bigg[ {n \atop k} \Bigg] = \frac{n!}{k!} \sum_{{i_1+i_2+...+i_k = n} \atop
	{0<i_1,...,i_k\leq n}} \frac{1}{i_1 i_2 ...i_k}
$$

\vspace{0.2cm}

\noindent \textit{Exercise 2.A.}
\textbf{Ad Cobweb tiling problem} [10] see also [11,12].
$\Big\{{\eta \atop \kappa}\Big\}_\lambda = ?$ i.e. the number of $\kappa$-block partitions
with block sizes all equal to $\lambda$ = ?
The answer is known  due to [7,8,9,10]. Prove it.

$$
	\Bigg\{ {\eta \atop \kappa} \Bigg\}_\lambda = \delta_{\eta,\kappa\lambda}
	\frac{\eta!}{\kappa!} \sum_{{i_1+i_2+...+i_\kappa = \eta} \atop {0<i_1=...
	=i_\kappa=\lambda}} \frac{1}{i_1 i_2 ...i_\kappa} =
	\delta_{\eta,\kappa\lambda} \frac{\eta!}{\kappa!(\lambda!)^\kappa}
$$

\noindent \textbf{Hint} to Exercise 1A.
$$|X^S| = |X|^|S| = k^n$$

\noindent therefore ($k$ summands)
$$k^n = (1+1+...+1)n = \sum_{i_1+i_2+...
+i_k=n}{n \choose {i_1,i_2,...,i_k}}$$

\noindent where $0 < i_1, i_2, ...,i_k = n$. A map from $S$
to $X$ is not a surjection  iff     $ i_1 = 0 \vee i_2=0 \vee ... \vee i_k = 0$.

\vspace{0.2cm}

\noindent \textbf{Consider Compositions} of the natural number $n$ via \textbf{partitions summation}

\vspace{0.2cm}

\noindent \textit{Exercise 3.A.}
Prove that
$$
	C(n,k,1,\infty) = {{n-1} \choose {k-1}} =
	\sum_{{1\lambda_1 +2\lambda_2 +...+n\lambda_n = n} \atop {\lambda_1+\lambda_2+...
	+\lambda_n = k}}
	\frac{k!}{\lambda_1!\lambda_2!...\lambda_n!}
$$

\section{Appendix : ``what is further on''}

\noindent \textbf{7.1.}
Recall.  The idea of graph compositions as introduced by A.\,Knopfmacher and M.\,E.\,Mays [2],
generalizes both ordinary compositions of positive integers and partitions of finite sets. In
[2--5] the authors provided various formulas, generating functions, and recurrence relations
for composition counting functions for several families of graphs.  Recently in [13]  some
of the results involving compositions of bipartite graphs  have been  derived in a simpler
way  using exponential generating functions.

\vspace{3mm}

\noindent \textbf{7.2.}
In the [14 ] a new construction of tree-like graphs where nodes are graphs themselves was
added and  examples of these tree-like compositions, a corresponding theorem and few
resulting conclusions are delivered. Compare with  [15].

\vspace{3mm}

\noindent \textbf{7.3.}
Ad naming and history roots  see and compare with [2] the content of references [16--19].
This comparison is important. For example, in [16]  W.H. Cunningham described a composition
for directed graphs, a composition which generalizes the substitution (or X-join) composition
of graphs and digraphs, as well as the graph version of set-family composition. There he s
proved that a general decomposition theory from [17] can be applied to the resulting digraph
decomposition. ``A consequence is a theorem which asserts the uniqueness of a decomposition
of any digraph, each member of the decomposition being either indecomposable or ``special''.
The special digraphs are completely characterized; they are members of a few interesting classes.
Efficient decomposition algorithms are also presented'' quoted from [16].
As for the Cunningham's split decomposition of an arbitrary undirected graph --
 the related composition operation is a generalization of the modular composition --
 also called substitution of X-join. The split decomposition is useful in recognizing special
 classes of graphs, such as circle graphs, which are the intersection graphs of arcs of a
 circle, and parity graphs, because these graphs are closed under the inverse composition
 operation. The decomposition can also be used to find NC algorithms for some optimization
 problems on special families of graphs, assuming these problems can be solved in NC for the
 indecomposable graphs of the decomposition.[see more in  \emph{http://hdl.handle.net/1813/6874}]
Let us illustrate the importance of   Cunningham`s split decomposition from [16] quoting as
an example the reference [20].  The authors of  [20]  introduce there a new  structural
property of parity graphs. Namely Cunningham`s split decomposition returns exactly, as
building blocks of parity graphs, cliques and bipartite graphs. This characterization and
the observation that the split decomposition process is performed in linear time, gives rise
to optimum algorithms for the recognition problem and for  the maximum weighted clique problem.

\vspace{0.5cm}

\noindent{\small Institute of Combinatorics and its Applications}

\noindent {\small High School of Mathematics and Applied Informatics}

\noindent {\small Kamienna 17, PL-15-021 Bia\l ystok}

\noindent {\small Poland}

\noindent {\small e-mail: kwandr@gmail.com}

\vspace{0.5cm}

\noindent Presented by Marek Moneta at the Session of the
Mathematical-Phys\-ical Commission of the \L \'od\'z Society of Sciences and
Arts on November 4, 2009

\vspace{0.8cm}

\noindent {\bf  O SK\L ADANIU LICZB I GRAF\'OW}

\vspace{0.2cm}

\noindent {\small S t r e s z c z e n i e}

{\small  G\l \'ownym celem obecnej noty jest zaproponowanie kilku
problem\'ow, kt\'ore daj\c{a} si\c{e} \l atwo sformu\l owa\'c, lecz
niekt\'ore z nich zapewne nie s\c{a} rozwi\c{a}zane. S\c{a} one o og\'olnym
znaczeniu dla matematyki dyskretnej w\l \c{a}czaj\c{a}c now\c{a}   ga\l
\c{a}zk\c{e} w ramach pewnej ga\l \c{e}zi teorii graf\'ow zwi\c{a}zanej z
ich sk\l adaniem. Wynikaj\c{a} st\c{a}d problemy, kt\'ore s\c{a}
istotne w zwi\c{a}zku z sze\-regiem zada\'n i wskaz\'owek opartych przede
wszystkim na ostatnich pracach Knopfmachera i wsp\'o\l autor\'ow.

}




\begin{thebibliography} {99}

\bibitem{1}  C.\,Kimberling,
\emph{Enumeration of paths, compositions of integers, and Fibonacci numbers},
The Fibonacci Quarterly, 39 (5), (2001), 430--435.

\bibitem{2} A.\,Knopfmacher and M.\,E.\,Mays,
\emph{Graph compositions I: Basic enumeration Integers},
The Electronic Journal of Combinatorial Number Theory  {\bf 1} (2001), A04.

\bibitem{3}  A.\,Knopfmacher and N.\,Robbins,
\emph{Compositions with parts constrained by the leading summand},
Ars Combinatoria {\bf 76} (2005), 287-295.

\bibitem{4}  A.\,Knopfmacher and M.\,E.\,Mays,
\emph{Compositions with m distinct parts},
Ars Combinatoria  {\bf 53}  (1999) 111--128.

\bibitem{5}  B.\,Richmont and  A.\,Knopfmacher,
\emph{Compositions with distinct parts},
Aequationes Math. {\bf 49} (1995), 86--97.

\bibitem{6}  T.\,Doslic,
\emph{Maximum Product Over Partitions Into Distinct Parts},
Journal of Integer Sequences  {\bf 8} (2005), Article 05.5.8.

\bibitem{7} Ch.\,Jordan,
\emph{On Stirling Numbers},
Tôhoku Math. J. {\bf 37}   (1933),  254--278,   formula (59).

\bibitem{8} C.\,Jordan,
\emph{Calculus of finite differences},
Chelsea, New York, 1947,   § 51, § 60.

\bibitem{9} T.\,Kreid,
\emph{Combinatorial sequences of polynomials},
Comment. Math. Prace Math. {\bf 29} (1990), 233--242.

\bibitem{10} A.\,K.\,Kwa\'sniewski,
\emph{Cobweb posets as noncommutative prefabs},
Adv. Stud. Contemp. Math.  {\bf 14} (1) (2007), 37--47.

\bibitem{11}  A.\,K.\,Kwa\'sniewski, \emph{On cobweb posets and their
combinatorially admissible sequences}, Adv. Studies Contemp. Math. Vol.
{\bf 18}, no.\,1 (2009), 17--32;  ArXiv:0512578v4 [v5] Mon, 19 Jan 2009
21:47:32 GMT.


\bibitem{12} M.\,Dziemia\'nczuk,
\emph{On Cobweb posets tiling problem},
arXiv: math.CO/0709.4263  4 Oct 2007.

\bibitem{13} A.\,Huq,
\emph{Compositions of Graphs Revisited},
The Electronic Journal of Combinatorics  {\bf 14} (2007), N15;
arXiv:math.CO/0704.3821, 28 Apr 2007.

\bibitem{14} W.\,Bajguz,
\emph{Graph and Union of Graphs Compositions}, Advanced Studies in Contemporary Math.
Vol. 16, No. 2, April 2008,  pp. 245-249.
arXiv:math.CO/0601755, 31 Jan 2006.

\bibitem{15} J.\,N.\,Ridley  and M.\,E.\,Mays,
\emph{Compositions  of Union of Graphs},
The Fibonacci Quarterly {\bf 42}, no.\,3 (2004), 222--230.

\bibitem{16} W.\,H.\,Cunningham,
\emph{Decomposition of Directed Graphs},
SIAM J. Alg. Disc. Meth. {\bf 3} (1982), 214--228.

\bibitem{17} W.\,H.\,Cunningham and J.\,Edmonds,
\emph{A Combinatorial Decomposition Theory},
Canad. J. Math. {\bf 32} (1980), 734--765.

\bibitem{18} L.\,O.\,James,  R.\,G.\,Stanton, and  D.\,D.\,Cowan,
\emph{Graph decomposition for undirected graphs},
Proceedings of the third Southeastern International Conference on Combinatorics,
Graph Theory, and Computing (CGTC'72) (1972), 281--290.

\bibitem{19} G.\,Sabidussi,
\emph{The composition of graphs},
Duke Math. J. {\bf 26} (1959), 693--696.

\bibitem{20} S.\,Cicerone and G.\,Di Stefano,
\emph{On the extension of bipartite to parity graphs},
Discrete Applied Mathematics  {\bf 95}, Issues 1--3, (1999), 181--195.


\end{thebibliography}
\end{document}